\documentclass[twocolumn]{article}

\usepackage{amsmath,amsthm,amssymb,amsfonts}
\usepackage{verbatim}

\usepackage{stmaryrd} 
\usepackage{hyperref}

\usepackage[colorinlistoftodos]{todonotes}

\usepackage{tikz}
\definecolor{steelblue}{RGB}{70,130,180}
\tikzset{vertex/.style={circle,draw,fill=steelblue,inner sep=0pt,minimum size=1mm}}
\tikzset{edge/.style={color=steelblue,line width=.5mm}}

\theoremstyle{plain}

\theoremstyle{definition}

\numberwithin{thm}{section}

\newcommand{\adjeq}{\leftrightarroweq}

\def\R{{\mathbb R}}

\title{Nice Neighbors: A brief adventure in mathematical gamification}
\author{P.~Christopher Staecker\thanks{Department of Mathematics, Fairfield University, Fairfield CT 06824, USA, \url{cstaecker@fairfield.edu}}}

\begin{document}
\bibliographystyle{plain}

\maketitle 
Last year I came across a strange graph theory problem in digital topology and decided to turn it into a little video game, to help wrap my mind around it. It turned out to be pretty fun, so I made it into a web game that other people could play. 
I took 3500 specific math problems and made each one into a level of the game, and I waited to see if people would play the game and solve my problems for me.
Within 2 months a lot of people and at least one non-person did play the game, and they did solve my problems for me. I'll try to describe the mathematics behind it as well as some of the surprises along the way that still have me scratching my head a bit. 

\begin{figure}[h]
\[ \includegraphics[scale=.1]{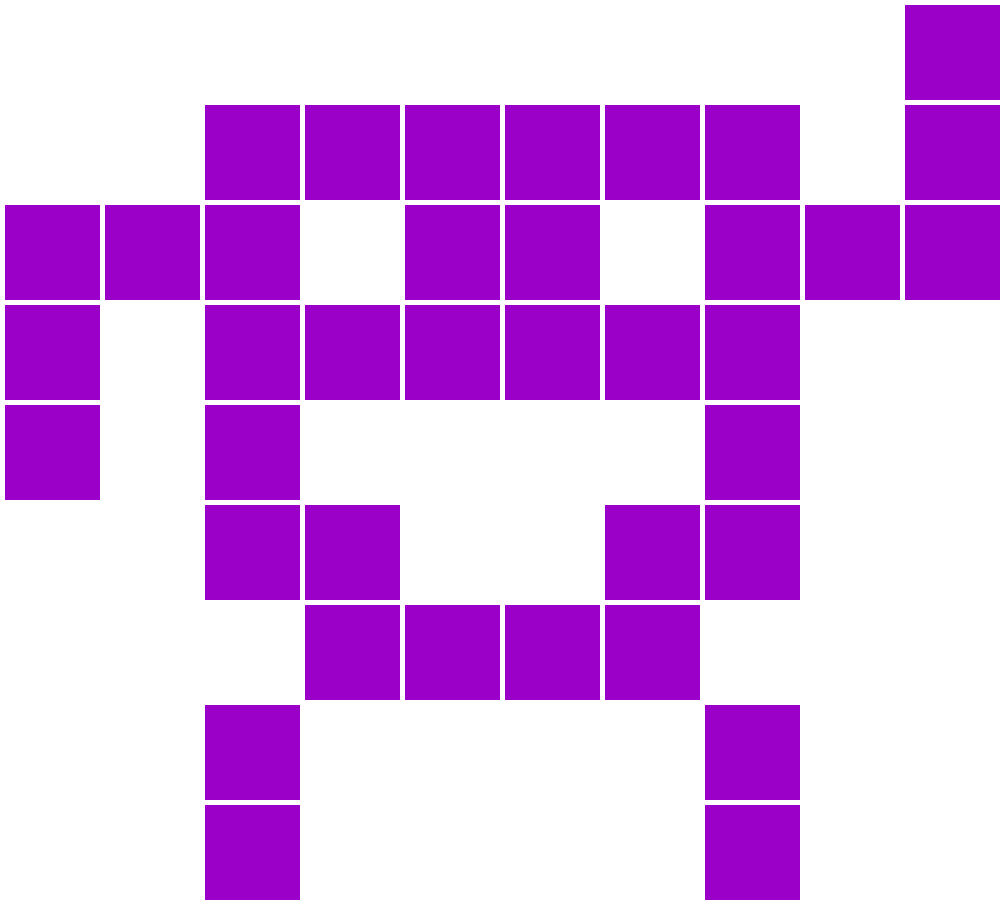} \]
\caption{Pappy, a friendly digital image. Fabulous pixel art by the author.\label{invaderfig}}
\end{figure}

\subsection*{A little digital topology}
My game is based on a mathematical question from the fairly new field of \emph{digital topology}, which is about topology for spaces made up of discrete pixels. Consider Pappy, the digital character in Figure \ref{invaderfig}, which we will think of as a set of 38 points. (Pappy was named with help from my young daughters, because he's purple and happy.)

Digital images like this don't exactly obey the traditional laws of Euclidean geometry and topology. For instance we can have two points of Pappy which are right next to each other with no point in between, which is impossible for points in $\R^2$. 


Today is an interesting time for digital topology. The mathematical field is still in its infancy, despite the fact that digital image analysis has been going on for decades in industry. It's a theme in the history of mathematics: useful ideas tend to appear in the ``real world'' before we mathematicians describe and define their theoretical foundations. Math historian Judith Grabiner said this about calculus, ``First the derivative was used, then discovered, explored and developed, and only then, defined.''

Today there are a few different approaches to digital toplogy which are actively being pursued. We're going to focus on a digital model that closely resembles graph theory, which seems to have its origins in the work of Rosenfeld in the 1970s. 
The basic idea is that the topological information in a digital image is entirely expressed by which pixels are adjacent to each other. 
Since all we care about are the pixels and their connections, we may as well think of a digital image as a graph with a vertex for each pixel, and an edge for each adjacency. 

\begin{figure}
\[ 
\begin{tikzpicture}[scale=.35]
\node[vertex,minimum size=1mm] (9) at (9,8) {};\node[vertex,minimum size=1mm] (12) at (2,7) {};\node[vertex,minimum size=1mm] (13) at (3,7) {};\node[vertex,minimum size=1mm] (14) at (4,7) {};\node[vertex,minimum size=1mm] (15) at (5,7) {};\node[vertex,minimum size=1mm] (16) at (6,7) {};\node[vertex,minimum size=1mm] (17) at (7,7) {};\node[vertex,minimum size=1mm] (19) at (9,7) {};\node[vertex,minimum size=1mm] (20) at (0,6) {};\node[vertex,minimum size=1mm] (21) at (1,6) {};\node[vertex,minimum size=1mm] (22) at (2,6) {};\node[vertex,minimum size=1mm] (24) at (4,6) {};\node[vertex,minimum size=1mm] (25) at (5,6) {};\node[vertex,minimum size=1mm] (27) at (7,6) {};\node[vertex,minimum size=1mm] (28) at (8,6) {};\node[vertex,minimum size=1mm] (29) at (9,6) {};\node[vertex,minimum size=1mm] (30) at (0,5) {};\node[vertex,minimum size=1mm] (32) at (2,5) {};\node[vertex,minimum size=1mm] (33) at (3,5) {};\node[vertex,minimum size=1mm] (34) at (4,5) {};\node[vertex,minimum size=1mm] (35) at (5,5) {};\node[vertex,minimum size=1mm] (36) at (6,5) {};\node[vertex,minimum size=1mm] (37) at (7,5) {};\node[vertex,minimum size=1mm] (40) at (0,4) {};\node[vertex,minimum size=1mm] (42) at (2,4) {};\node[vertex,minimum size=1mm] (47) at (7,4) {};\node[vertex,minimum size=1mm] (52) at (2,3) {};\node[vertex,minimum size=1mm] (53) at (3,3) {};\node[vertex,minimum size=1mm] (56) at (6,3) {};\node[vertex,minimum size=1mm] (57) at (7,3) {};\node[vertex,minimum size=1mm] (63) at (3,2) {};\node[vertex,minimum size=1mm] (64) at (4,2) {};\node[vertex,minimum size=1mm] (65) at (5,2) {};\node[vertex,minimum size=1mm] (66) at (6,2) {};\node[vertex,minimum size=1mm] (72) at (2,1) {};\node[vertex,minimum size=1mm] (77) at (7,1) {};\node[vertex,minimum size=1mm] (82) at (2,0) {};\node[vertex,minimum size=1mm] (87) at (7,0) {};\draw[edge,line width=.5mm] (12) -- (13);\draw[edge,line width=.5mm] (13) -- (14);\draw[edge,line width=.5mm] (14) -- (15);\draw[edge,line width=.5mm] (15) -- (16);\draw[edge,line width=.5mm] (16) -- (17);\draw[edge,line width=.5mm] (19) -- (9);\draw[edge,line width=.5mm] (20) -- (21);\draw[edge,line width=.5mm] (21) -- (22);\draw[edge,line width=.5mm] (22) -- (12);\draw[edge,line width=.5mm] (24) -- (25);\draw[edge,line width=.5mm] (24) -- (14);\draw[edge,line width=.5mm] (25) -- (15);\draw[edge,line width=.5mm] (27) -- (28);\draw[edge,line width=.5mm] (27) -- (17);\draw[edge,line width=.5mm] (28) -- (29);\draw[edge,line width=.5mm] (29) -- (19);\draw[edge,line width=.5mm] (30) -- (20);\draw[edge,line width=.5mm] (32) -- (33);\draw[edge,line width=.5mm] (32) -- (22);\draw[edge,line width=.5mm] (33) -- (34);\draw[edge,line width=.5mm] (34) -- (35);\draw[edge,line width=.5mm] (34) -- (24);\draw[edge,line width=.5mm] (35) -- (36);\draw[edge,line width=.5mm] (35) -- (25);\draw[edge,line width=.5mm] (36) -- (37);\draw[edge,line width=.5mm] (37) -- (27);\draw[edge,line width=.5mm] (40) -- (30);\draw[edge,line width=.5mm] (42) -- (32);\draw[edge,line width=.5mm] (47) -- (37);\draw[edge,line width=.5mm] (52) -- (53);\draw[edge,line width=.5mm] (52) -- (42);\draw[edge,line width=.5mm] (56) -- (57);\draw[edge,line width=.5mm] (57) -- (47);\draw[edge,line width=.5mm] (63) -- (64);\draw[edge,line width=.5mm] (63) -- (53);\draw[edge,line width=.5mm] (64) -- (65);\draw[edge,line width=.5mm] (65) -- (66);\draw[edge,line width=.5mm] (66) -- (56);\draw[edge,line width=.5mm] (82) -- (72);\draw[edge,line width=.5mm] (87) -- (77);
\end{tikzpicture}
\qquad
\begin{tikzpicture}[scale=.35]
\node[vertex,minimum size=1mm] (9) at (9,8) {};\node[vertex,minimum size=1mm] (12) at (2,7) {};\node[vertex,minimum size=1mm] (13) at (3,7) {};\node[vertex,minimum size=1mm] (14) at (4,7) {};\node[vertex,minimum size=1mm] (15) at (5,7) {};\node[vertex,minimum size=1mm] (16) at (6,7) {};\node[vertex,minimum size=1mm] (17) at (7,7) {};\node[vertex,minimum size=1mm] (19) at (9,7) {};\node[vertex,minimum size=1mm] (20) at (0,6) {};\node[vertex,minimum size=1mm] (21) at (1,6) {};\node[vertex,minimum size=1mm] (22) at (2,6) {};\node[vertex,minimum size=1mm] (24) at (4,6) {};\node[vertex,minimum size=1mm] (25) at (5,6) {};\node[vertex,minimum size=1mm] (27) at (7,6) {};\node[vertex,minimum size=1mm] (28) at (8,6) {};\node[vertex,minimum size=1mm] (29) at (9,6) {};\node[vertex,minimum size=1mm] (30) at (0,5) {};\node[vertex,minimum size=1mm] (32) at (2,5) {};\node[vertex,minimum size=1mm] (33) at (3,5) {};\node[vertex,minimum size=1mm] (34) at (4,5) {};\node[vertex,minimum size=1mm] (35) at (5,5) {};\node[vertex,minimum size=1mm] (36) at (6,5) {};\node[vertex,minimum size=1mm] (37) at (7,5) {};\node[vertex,minimum size=1mm] (40) at (0,4) {};\node[vertex,minimum size=1mm] (42) at (2,4) {};\node[vertex,minimum size=1mm] (47) at (7,4) {};\node[vertex,minimum size=1mm] (52) at (2,3) {};\node[vertex,minimum size=1mm] (53) at (3,3) {};\node[vertex,minimum size=1mm] (56) at (6,3) {};\node[vertex,minimum size=1mm] (57) at (7,3) {};\node[vertex,minimum size=1mm] (63) at (3,2) {};\node[vertex,minimum size=1mm] (64) at (4,2) {};\node[vertex,minimum size=1mm] (65) at (5,2) {};\node[vertex,minimum size=1mm] (66) at (6,2) {};\node[vertex,minimum size=1mm] (72) at (2,1) {};\node[vertex,minimum size=1mm] (77) at (7,1) {};\node[vertex,minimum size=1mm] (82) at (2,0) {};\node[vertex,minimum size=1mm] (87) at (7,0) {};\draw[edge,line width=.5mm] (12) -- (13);\draw[edge,line width=.5mm] (13) -- (14);\draw[edge,line width=.5mm] (13) -- (24);\draw[edge,line width=.5mm] (14) -- (15);\draw[edge,line width=.5mm] (14) -- (25);\draw[edge,line width=.5mm] (15) -- (16);\draw[edge,line width=.5mm] (16) -- (17);\draw[edge,line width=.5mm] (16) -- (27);\draw[edge,line width=.5mm] (17) -- (28);\draw[edge,line width=.5mm] (19) -- (9);\draw[edge,line width=.5mm] (20) -- (21);\draw[edge,line width=.5mm] (21) -- (22);\draw[edge,line width=.5mm] (21) -- (12);\draw[edge,line width=.5mm] (21) -- (32);\draw[edge,line width=.5mm] (22) -- (12);\draw[edge,line width=.5mm] (22) -- (13);\draw[edge,line width=.5mm] (22) -- (33);\draw[edge,line width=.5mm] (24) -- (25);\draw[edge,line width=.5mm] (24) -- (14);\draw[edge,line width=.5mm] (24) -- (15);\draw[edge,line width=.5mm] (24) -- (35);\draw[edge,line width=.5mm] (25) -- (15);\draw[edge,line width=.5mm] (25) -- (16);\draw[edge,line width=.5mm] (25) -- (36);\draw[edge,line width=.5mm] (27) -- (28);\draw[edge,line width=.5mm] (27) -- (17);\draw[edge,line width=.5mm] (28) -- (29);\draw[edge,line width=.5mm] (28) -- (19);\draw[edge,line width=.5mm] (29) -- (19);\draw[edge,line width=.5mm] (30) -- (20);\draw[edge,line width=.5mm] (30) -- (21);\draw[edge,line width=.5mm] (32) -- (33);\draw[edge,line width=.5mm] (32) -- (22);\draw[edge,line width=.5mm] (33) -- (34);\draw[edge,line width=.5mm] (33) -- (24);\draw[edge,line width=.5mm] (34) -- (35);\draw[edge,line width=.5mm] (34) -- (24);\draw[edge,line width=.5mm] (34) -- (25);\draw[edge,line width=.5mm] (35) -- (36);\draw[edge,line width=.5mm] (35) -- (25);\draw[edge,line width=.5mm] (36) -- (37);\draw[edge,line width=.5mm] (36) -- (27);\draw[edge,line width=.5mm] (36) -- (47);\draw[edge,line width=.5mm] (37) -- (27);\draw[edge,line width=.5mm] (37) -- (28);\draw[edge,line width=.5mm] (40) -- (30);\draw[edge,line width=.5mm] (42) -- (32);\draw[edge,line width=.5mm] (42) -- (33);\draw[edge,line width=.5mm] (42) -- (53);\draw[edge,line width=.5mm] (47) -- (37);\draw[edge,line width=.5mm] (52) -- (53);\draw[edge,line width=.5mm] (52) -- (42);\draw[edge,line width=.5mm] (52) -- (63);\draw[edge,line width=.5mm] (53) -- (64);\draw[edge,line width=.5mm] (56) -- (57);\draw[edge,line width=.5mm] (56) -- (47);\draw[edge,line width=.5mm] (57) -- (47);\draw[edge,line width=.5mm] (63) -- (64);\draw[edge,line width=.5mm] (63) -- (53);\draw[edge,line width=.5mm] (64) -- (65);\draw[edge,line width=.5mm] (65) -- (66);\draw[edge,line width=.5mm] (65) -- (56);\draw[edge,line width=.5mm] (66) -- (56);\draw[edge,line width=.5mm] (66) -- (57);\draw[edge,line width=.5mm] (66) -- (77);\draw[edge,line width=.5mm] (72) -- (63);\draw[edge,line width=.5mm] (82) -- (72);\draw[edge,line width=.5mm] (87) -- (77);
\end{tikzpicture}
\]
\caption{Two graphs modeling Pappy. On the left using 4-adjacency, and on the right using 8-adjacency.\label{invader48fig}}
\end{figure}
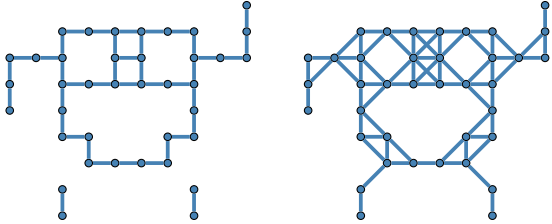

What exactly do we mean by ``adjacent''? In the plane there are two standard schemes: ``4-adjacency'' which allows no diagonals, or ``8-adjacency'' which does allow diagonals. 
Our friend Pappy, for example, would be represented by the graphs in Figure \ref{invader48fig}, which look like two versions of Pappy's skeleton. (Says one daughter: ``That's what Pappy will look like when he's dead.'' The other adds, ``He still looks happy.'')
As you can see, the topological structure of the image is pretty different depending on which adjacency we use. 
This is a bit of a contrast with the several standard metrics on $\R^2$ which all give the same topology.

\subsection*{Reducible images}
In the Summer of 2014 I was advising an REU project at Fairfield University with my students Jason Haarmann, Meg Fields, and Casey Peters. 
The idea we stumbled upon has to do with ``reducing'' a digital image by removing pixels which don't affect the topology. 
Just like in ordinary topology, we consider information about length or angle as not so important, while information about holes or connectivity is what we really care about. 
For example the pixels on the ends of Pappy's hands should be removable without changing the topology.

Here's the formal idea that we settled on: When two pixels $x$ and $y$ are adjacent or equal, let's write $x\adjeq y$. We say a digital image $X$ is \emph{reducible} if there is some function $f:X\to X$ such that:
\begin{itemize}
\item $f$ is not onto,
\item $x \adjeq f(x)$ for every pixel $x$,
\item $f(x)\adjeq f(y)$ whenever $x\adjeq y$.
\end{itemize}
(In digital topology jargon, $f$ gives a \emph{homotopy equivalence}.) If no such $f$ exists, then we call the image \emph{irreducible}. You should think of $f$ as rearranging the vertices of the graph. Then the three conditions say that: one vertex spot must be vacant after the move, each vertex can only move to an adjacent spot, and any adjacent vertices must still be adjacent after the move. 

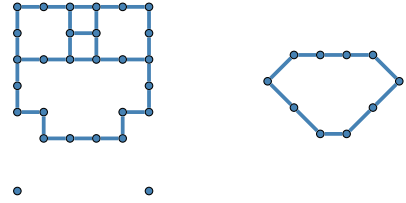
\begin{figure}
\[ 
\begin{tikzpicture}[scale=.35]
\node[vertex,minimum size=1mm] (12) at (2,7) {};\node[vertex,minimum size=1mm] (13) at (3,7) {};\node[vertex,minimum size=1mm] (14) at (4,7) {};\node[vertex,minimum size=1mm] (15) at (5,7) {};\node[vertex,minimum size=1mm] (16) at (6,7) {};\node[vertex,minimum size=1mm] (17) at (7,7) {};\node[vertex,minimum size=1mm] (22) at (2,6) {};\node[vertex,minimum size=1mm] (24) at (4,6) {};\node[vertex,minimum size=1mm] (25) at (5,6) {};\node[vertex,minimum size=1mm] (27) at (7,6) {};\node[vertex,minimum size=1mm] (32) at (2,5) {};\node[vertex,minimum size=1mm] (33) at (3,5) {};\node[vertex,minimum size=1mm] (34) at (4,5) {};\node[vertex,minimum size=1mm] (35) at (5,5) {};\node[vertex,minimum size=1mm] (36) at (6,5) {};\node[vertex,minimum size=1mm] (37) at (7,5) {};\node[vertex,minimum size=1mm] (42) at (2,4) {};\node[vertex,minimum size=1mm] (47) at (7,4) {};\node[vertex,minimum size=1mm] (52) at (2,3) {};\node[vertex,minimum size=1mm] (53) at (3,3) {};\node[vertex,minimum size=1mm] (56) at (6,3) {};\node[vertex,minimum size=1mm] (57) at (7,3) {};\node[vertex,minimum size=1mm] (63) at (3,2) {};\node[vertex,minimum size=1mm] (64) at (4,2) {};\node[vertex,minimum size=1mm] (65) at (5,2) {};\node[vertex,minimum size=1mm] (66) at (6,2) {};\node[vertex,minimum size=1mm] (82) at (2,0) {};\node[vertex,minimum size=1mm] (87) at (7,0) {};\draw[edge,line width=.5mm] (12) -- (13);\draw[edge,line width=.5mm] (13) -- (14);\draw[edge,line width=.5mm] (14) -- (15);\draw[edge,line width=.5mm] (15) -- (16);\draw[edge,line width=.5mm] (16) -- (17);\draw[edge,line width=.5mm] (22) -- (12);\draw[edge,line width=.5mm] (24) -- (25);\draw[edge,line width=.5mm] (24) -- (14);\draw[edge,line width=.5mm] (25) -- (15);\draw[edge,line width=.5mm] (27) -- (17);\draw[edge,line width=.5mm] (32) -- (33);\draw[edge,line width=.5mm] (32) -- (22);\draw[edge,line width=.5mm] (33) -- (34);\draw[edge,line width=.5mm] (34) -- (35);\draw[edge,line width=.5mm] (34) -- (24);\draw[edge,line width=.5mm] (35) -- (36);\draw[edge,line width=.5mm] (35) -- (25);\draw[edge,line width=.5mm] (36) -- (37);\draw[edge,line width=.5mm] (37) -- (27);\draw[edge,line width=.5mm] (42) -- (32);\draw[edge,line width=.5mm] (47) -- (37);\draw[edge,line width=.5mm] (52) -- (53);\draw[edge,line width=.5mm] (52) -- (42);\draw[edge,line width=.5mm] (56) -- (57);\draw[edge,line width=.5mm] (57) -- (47);\draw[edge,line width=.5mm] (63) -- (64);\draw[edge,line width=.5mm] (63) -- (53);\draw[edge,line width=.5mm] (64) -- (65);\draw[edge,line width=.5mm] (65) -- (66);\draw[edge,line width=.5mm] (66) -- (56);
\end{tikzpicture}
\qquad
\begin{tikzpicture}[scale=.35]
\node[vertex,minimum size=1mm] (24) at (4,5) {};\node[vertex,minimum size=1mm] (25) at (5,5) {};\node[vertex,minimum size=1mm] (33) at (3,5) {};\node[vertex,minimum size=1mm] (36) at (6,5) {};\node[vertex,minimum size=1mm] (42) at (2,4) {};\node[vertex,minimum size=1mm] (47) at (7,4) {};\node[vertex,minimum size=1mm] (53) at (3,3) {};\node[vertex,minimum size=1mm] (56) at (6,3) {};\node[vertex,minimum size=1mm] (64) at (4,2) {};\node[vertex,minimum size=1mm] (65) at (5,2) {};\draw[edge,line width=.5mm] (24) -- (25);\draw[edge,line width=.5mm] (25) -- (36);\draw[edge,line width=.5mm] (33) -- (24);\draw[edge,line width=.5mm] (36) -- (47);\draw[edge,line width=.5mm] (42) -- (33);\draw[edge,line width=.5mm] (42) -- (53);\draw[edge,line width=.5mm] (53) -- (64);\draw[edge,line width=.5mm] (56) -- (47);\draw[edge,line width=.5mm] (64) -- (65);\draw[edge,line width=.5mm] (65) -- (56);
	\node[draw=none] () at (7,7) {};
	\node[draw=none] () at (0,0) {};
\end{tikzpicture}
\]
\caption{Minimal reductions of Pappy using either 4- or 8-adjacency.\label{pappyreduced}}
\end{figure}

Our little buddy Pappy is reducible in several ways: we can chop off the ends of his arms or feet, for example. If we keep on reducing Pappy as much as possible, we obtain the irreducible graphs in Figure \ref{pappyreduced}.

\begin{figure}
\begin{align*}
\begin{tikzpicture}[scale=.7]
	\node[vertex] ()   at (0,0) {};
	\node[draw=none] () at (1,1) {};
	\node[draw=none] () at (-1,-1) {};
\end{tikzpicture}
\qquad
\begin{tikzpicture}[scale=.7]
	\node[vertex] (a)   at (0,1) {};
	\node[vertex] (b)   at (162:1) {}; 
	\node[vertex] (c) at (234:1) {};
	\node[vertex] (d) at (306:1)   {};
	\node[vertex] (e) at (18:1)  {};		
	\draw[edge] (a) -- (b) -- (c) -- (d) -- (e) -- (a);
	\node[draw=none] () at (1,1) {};
	\node[draw=none] () at (-1,-1) {};
\end{tikzpicture}
\qquad
\begin{tikzpicture}[scale=.7]
	\node[vertex](a) at (0,1) {};
	\node[vertex](b) at (150:1) {};
	\node[vertex](c) at (210:1) {};
	\node[vertex](d) at (270:1) {};
	\node[vertex](e) at (330:1) {};
	\node[vertex](f) at (30:1) {};
	\draw[edge] (a) -- (b) -- (c) -- (d) -- (e) -- (f) -- (a);
	\node[draw=none] () at (1,1) {};
	\node[draw=none] () at (-1,-1) {};
\end{tikzpicture}
\\
\begin{tikzpicture}[scale=.7]
	\node[vertex](a) at (0,1) {};
	\node[vertex](b) at (141:1) {};
	\node[vertex](c) at (193:1) {};
	\node[vertex](d) at (244:1) {};
	\node[vertex](e) at (296:1) {};
	\node[vertex](f) at (347:1) {};
	\node[vertex](g) at (39:1) {};
	\draw[edge] (a) -- (b) -- (c) -- (d) -- (e) -- (f) -- (g) -- (a);
	\node[draw=none] () at (1,1) {};
	\node[draw=none] () at (-1,-1) {};
\end{tikzpicture}
\qquad
\begin{tikzpicture}[scale=.7]
	\node[vertex](a) at (0,1) {};
	\node[vertex](b) at (150:1) {};
	\node[vertex](c) at (210:1) {};
	\node[vertex](d) at (270:1) {};
	\node[vertex](e) at (330:1) {};
	\node[vertex](f) at (30:1) {};
	\node[vertex](g) at (0,0) {};
	\draw[edge] (a) -- (b) -- (c) -- (d) -- (e) -- (f) -- (a);
	\draw[edge] (a) -- (g) -- (d);
	\node[draw=none] () at (1,1) {};
	\node[draw=none] () at (-1,-1) {};
\end{tikzpicture}
\qquad
\begin{tikzpicture}[scale=.7]
	\node[vertex] (a)   at (0,-1) {};
	\node[vertex] (b)   at (-162:1) {}; 
	\node[vertex] (c) at (-234:1) {};
	\node[vertex] (d) at (-306:1)   {};
	\node[vertex] (e) at (-18:1)  {};
	\node[vertex] (f) at (0,-.33) {};
	\node[vertex] (g) at (0,.33) {};		
	\draw[edge] (a) -- (b) -- (c) -- (d) -- (e) -- (a);
	\draw[edge] (d) -- (g) -- (c);
	\draw[edge] (g) -- (f) -- (a);
	\node[draw=none] () at (1,1) {};
	\node[draw=none] () at (-1,-1) {};
\end{tikzpicture}
\end{align*}
\caption{Catalog of all irreducible graphs up to 7 points\label{catalog}.}
\end{figure}
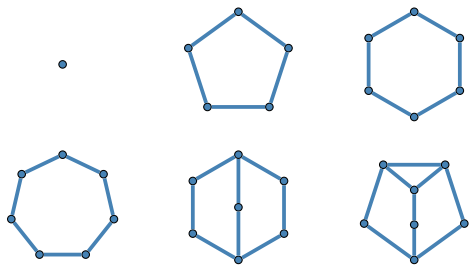

One major question that we tackled in our REU project was: which graphs are reducible, and which are not? Lacking any really big ideas, we started small: we managed to determine the complete list of connected irreducible graphs having at most 7 points. The result is Figure \ref{catalog}. This involved proving that these particular graphs are irreducible, and that \emph{all} other graphs of 7 points or less are reducible. (Note that we are considering any graphs at all, not just ones which come from 4- or 8-adjacency digital images in the plane. The question: ``given a graph, decide efficiently whether or not it is realizable as the graph of a 4- or 8-adjacency digital image in the plane'' is itself an open problem that we didn't want to get into.)

We stopped at 7 points because things got too complicated: we had to use a computer search to rule out thousands of graphs, and then had to check 15 special cases by hand that the computer couldn't handle. For 8 points the number of special cases was 160, and for 9 points it was 3251. 

\subsection*{The game}
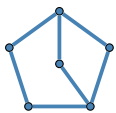
\begin{figure}
\[
\begin{tikzpicture}[scale=.7]
	\node[vertex] (a)   at (0,1) {};
	\node[vertex] (b)   at (162:1) {}; 
	\node[vertex] (c) at (234:1) {};
	\node[vertex] (d) at (306:1)   {};
	\node[vertex] (e) at (18:1)  {};
	\node[vertex] (f) at (0,0) {};		
	\draw[edge] (a) -- (b) -- (c) -- (d) -- (e) -- (a);
	\draw[edge] (a) -- (f) -- (d);
\end{tikzpicture}
\]
\caption{There is a reduction for this graph, but it had our research group stumped for a while.\label{twistfig}}
\end{figure}

Deciding whether or not a given graph is reducible can be hard-- try Figure \ref{twistfig}. I tried to visualize the problem like one of those little board games where you have pegs that can only move between adjacent holes. But I couldn't ever think of a physical representation that matched the three rules for the reducing function $f$. Plus I'm terrible at those peg games.

A few months later, it occured to me in the shower (where I get most of my good mathematical ideas) that this could work as a computer game. I imagined a game which keeps track of the original graph, but allows you to move the vertices around like pegs in holes. The game would only let you move pegs in certain ways to make sure your rearrangement satisfies the 3 properties.

\begin{figure*}
\begin{align*}
\includegraphics[width=150px]{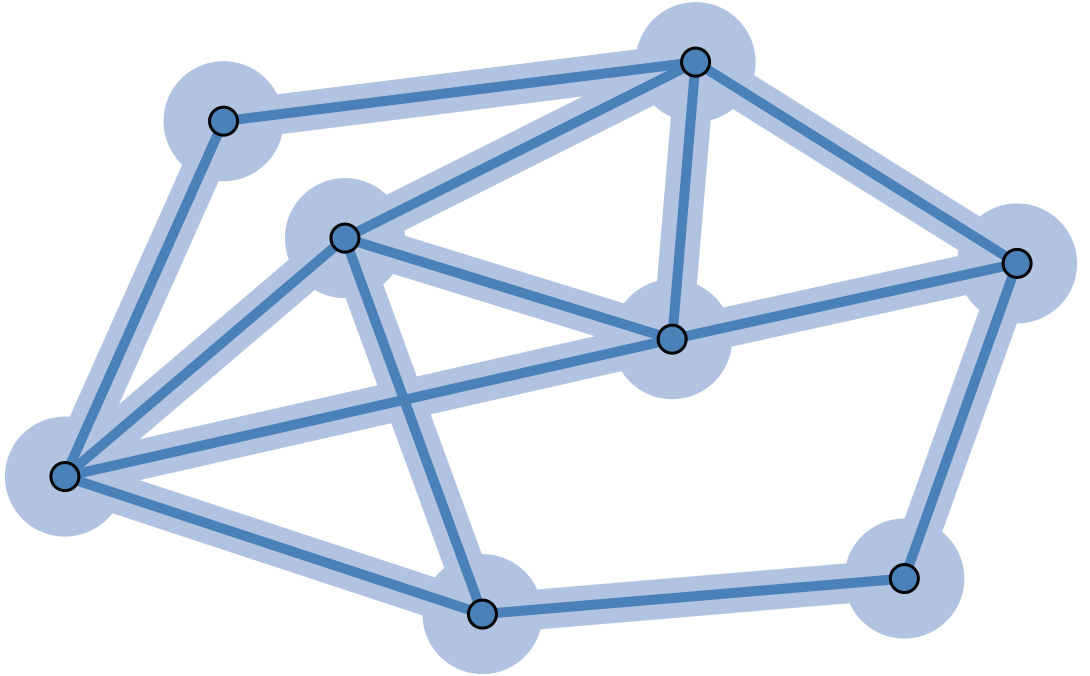}
\quad 
\includegraphics[width=150px]{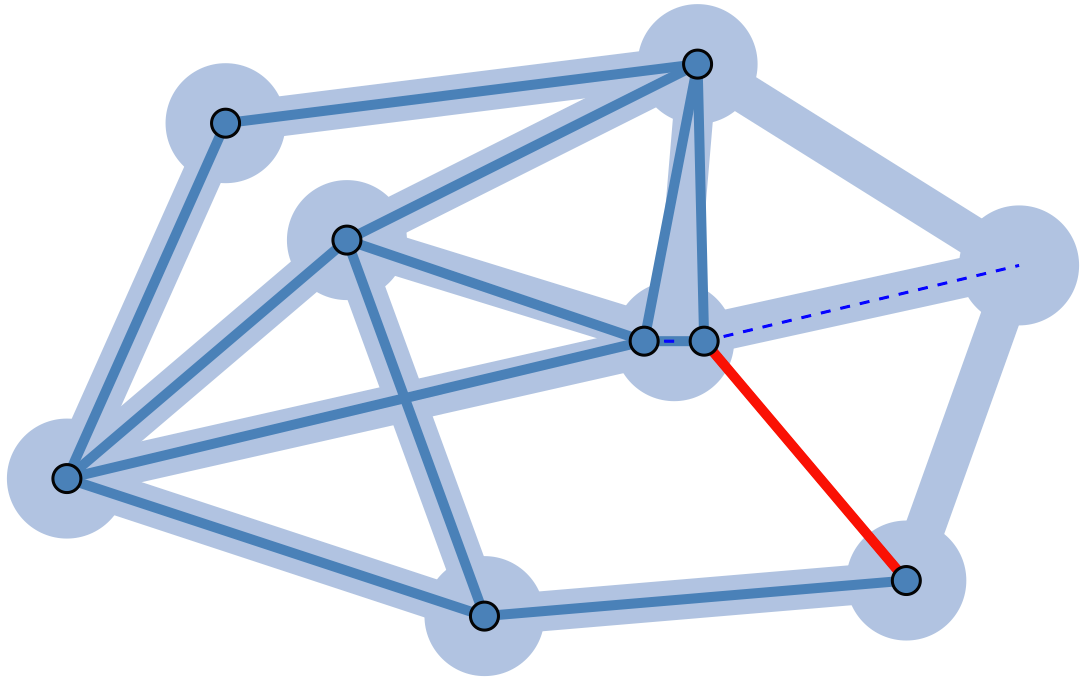}
\quad
\includegraphics[width=150px]{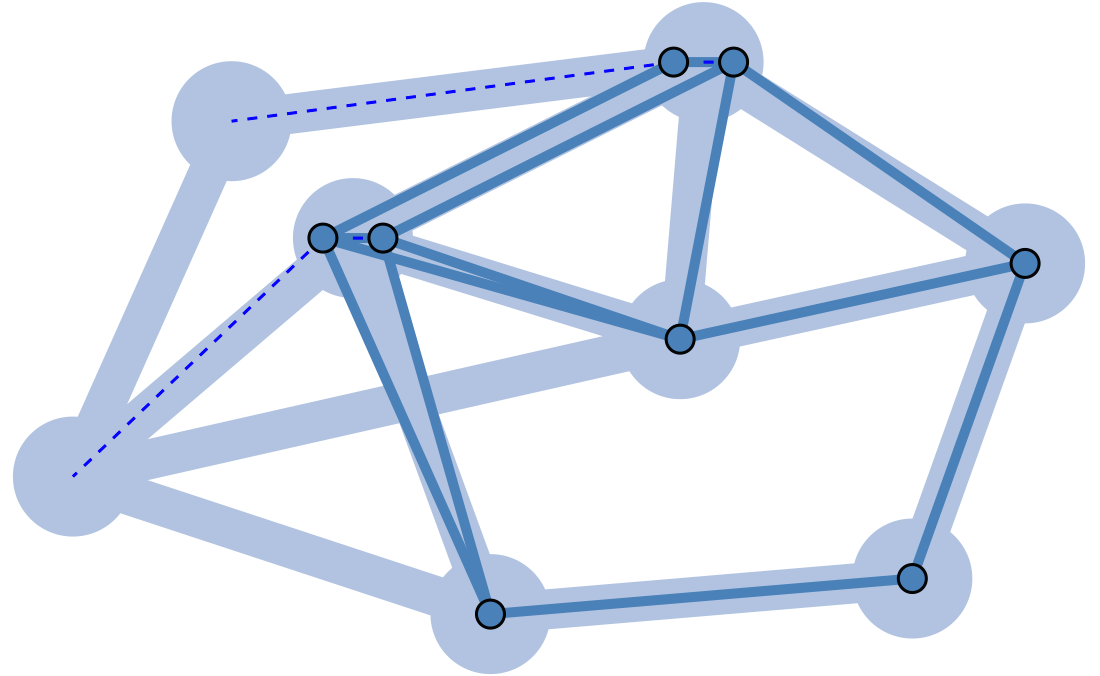}
\end{align*}
\caption{At left, a typical game level before any moves are made. At the middle, the player has dragged one vertex across to another spot. The red edge indicates that two vertices which should be adjacent are no longer adjacent, so this is not a winning position. The dotted blue line reminds the player where each vertex started. At right is a winning position for this level.\label{screenshots}}
\end{figure*}

I'm a decent programmer, and winter break was just starting, so I got to work. A few weeks later I had a working web game which gave a satisfying way to physically manipulate these graphs to get a feel for them. Figure \ref{screenshots} shows a typical game level.

In my game you are given a graph, and you win by dragging the vertices around onto one another so that:
\begin{itemize}
\item one spot is vacant,
\item each vertex ends up adjacent to where it started,
\item any vertices adjacent at the beginning are still adjacent at the end.
\end{itemize}
These three conditions correspond exactly to our 3 conditions for a reduction, so winning the game constitutes a proof that the graph is reducible.

I was a bit surprised to find that the game was pretty fun, and that it was possible to get better with practice. This meant that I was developing new intuition about my math problem, though it was usually hard to articulate exactly what I was learning.

Around this time it occured to me that I could code in all our thousands of special cases for 8 and 9 points and crowd-source the whole thing. I'd originally thought of the game as just being for myself, but it seemed natural to open it up. So I came up with the deeply unsatisfying name ``Nice Neighbors'', and cleaned up the interface a bit. 


Marketing the game was easier than I'd expected. 
I managed to convince some high-profile math tweeters to mention it, and within a couple of weeks I had been written up on some blogs and I had all the traffic I needed.

Pretty quickly I was getting feedback about the game, either directly from people who emailed me, or from reading tweets and blog comments. I learned that most people won't read the instructions, and also that people will figure out how to cheat (especially when the game code is open). I also saw that most of the levels were solved by a small group of ``hard-core'' Nice Neighbors enthusiasts. 





The biggest surprise of all was a particular player who I'll call User 87. (Full name: User 8709884216. I gave all users a random 10-digit identifier so I could track how many were coming from each person, and when.) Over the course of a day and a half, User 87 solved about one level per minute, until there were no more levels to solve. One level per minute, especially on the 9-point graphs, is quite a bit faster than I could consistently beat these levels, and User 87 wasn't stopping to sleep.

The only explanation I have for this is that User 87 was a computer algorithm which had bypassed the game's front-end to submit the solutions directly to my database.
Had I been trolled? Was User 87 trying to humiliate me? Or maybe just trying to help? And why didn't they ever contact me? I still don't know. Just to be safe, I wrote my own algorithms to reproduce User 87's solutions and make sure they were all correct (they were). 



In any event, within about 2 months from the game's launch, I had a complete catalog of irreducible graphs up through 9 points. 
I was even able to fulfill one of my career dreams-- some sequences in the great Online Encyclopedia of Integer Sequences. OEIS sequence A248571 is the number of irreducible graphs on $n$ points (starting with $n=1$): 1,0,0,0,1,1,3,28,547. The last two terms come from the game.
But how should I feel about all this? The takeaway is: if you ask the internet to solve a few thousand math puzzles, it probably will. But did I really learn anything, or gain new insight into my problem? Not exactly-- I just have all the answers now. 

Good mathematics isn't just about finding answers, but about explaining things. There's a big difference between answers and explanation-- these crowd-provided answers don't really illuminate anything. It's satisfying for me to have the answers, but I can't escape the feeling that I still don't understand why these answers are true. 

Hopefully somebody will figure it all out some day. Maybe you? Maybe me? Or maybe I'll just play my fun little game some more\dots

\subsection*{References and further reading}
A good place to start reading about Rosenfeld-style digital topology is Rosenfeld's original paper \cite{rose79}, or some later important papers by Boxer, e.g. \cite{boxe94}. Another important (and completely different) model of digital topology is based on the work of Khalimsky, which is covered nicely in the textbook by Adams and Franzosa \cite{af07}. Our REU team published a paper about reducibility and homotopy equivalence \cite{hmps15}, and you can read more about my algorithms in \cite{stae15}. 

Nice Neighbors is still up and running for your enjoyment \cite{nn}, as is the great Planarity \cite{planarity} which was my major design inspiration.

\renewcommand{\section}[2]{} 

\subsection*{Author bio}
Chris Staecker is Associate Professor in the Mathematics Department at Fairfield University, where he also teaches some computer programming. He enjoys boring movies and mechanical calculating machines.

\subsection*{Reducing Figure \ref{twistfig}}
To reduce the graph in Figure \ref{twistfig}, move each of the outer 5 points ``counter clockwise'' by 1 position, and move the center point to the top. We proved in the REU paper that any reduction of this graph must move all of its points, which answered an open question posed by Boxer.

\end{document}